\documentclass[12pt]{article}

\usepackage{amssymb,amsthm,amsfonts,graphicx,color,enumerate,float,subfigure,comment}

\newtheorem{thm}{Theorem}%[section]
\newtheorem{lemma}[thm]{Lemma}
\newtheorem{cor}[thm]{Corollary}
\newtheorem{prop}[thm]{Proposition}

\theoremstyle{remark}
\newtheorem{remark}[thm]{Remark}
\newtheorem{important remark}[thm]{Important Remark}

\def\R{{\mathbb R}}

\def\F{{\mathbb F}}

\title{Quasi-homomorphisms on mapping class groups}
  \author{Mladen Bestvina and Koji Fujiwara 
\thanks{The first author
  gratefully acknowledges the support by the National Science
  Foundation. }}

\begin{document}

\maketitle 

\centerline{\it To Sibe Marde\v si\' c, with admiration}

\abstract{We refine the construction of quasi-homomorphisms
on mapping class groups. It is useful to know that there are unbounded
quasi-homomorphisms which are bounded when restricted to particular
subgroups since then one deduces that the mapping class group is not
boundedly generated by these subgroups. In this note we enlarge the
class of such subgroups. The generalization is motivated by
considerations in first order theory of free groups.}

\section{Introduction}

Recall that a
{\it quasi-homomorphism} on a group $G$ is a function $h:G\to \R$ such
that
$$\sup_{\gamma_1,\gamma_2\in
  G}|h(\gamma_1\gamma_2)-h(\gamma_1)-h(\gamma_2)|<\infty$$
The set of all quasi-homomorphisms on $G$ forms a vector space
  ${\mathcal V}(G)$. Any bounded function on $G$ is trivially a
  quasi-homomorphism, and we consider the vector space $QH(G)$
  which is the quotient of ${\mathcal V}(G)$ by the subspace of
  bounded functions. The existence of many quasi-homomorphisms
  $G\to\R$ has implications for (stable) commutator length
  \cite{bavard}, (failure of) bounded generation of $G$, and
  non-embedding results for arithmetic lattices. For example, if
  $h:G\to\R$ is an unbounded quasi-homomorphism such that the
  restrictions $h|Q_i$ to subgroups $Q_1,\cdots,Q_k$ are bounded, then
  $G$ is not boundedly generated by the $Q_i$'s, i.e. for no $N>0$ can
  every element of $G$ be written as the product of $\leq N$ elements
  of $Q_1\cup\cdots\cup Q_k$. 

In this paper we consider the case when $G$ is the mapping class group
$MCG(S)$ 
of a compact surface $S$ (this is the group of isotopy classes of
homeomorphisms $S\to S$). With the exception of a small number of sporadic
surfaces, $MCG(S)$ is a ``large'' group, e.g. it contains a nonabelian
free subgroup.
In \cite{bf} we proved that $QH(MCG(S))$ is infinite-dimensional
provided $MCG(S)$ is not virtually abelian (and we only considered the
case of orientable $S$). More precisely, we have

\begin{thm}\label{old}
Let $S$ be a compact surface such that $MCG(S)$ is not virtually
abelian. Then $QH(MCG(S))$ is infinite-dimensional. Moreover, for
every finite collection of cyclic subgroups $$C_1,\cdots,C_k\subset
MCG(S)$$ there is an infinite-dimensional subspace of $QH(MCG(S))$ such
that each representative quasi-homomorphism $h$ of any element in this
subspace satisfies:
\begin{enumerate}[(1)]
\item $h$ is bounded on each $C_i$, and
\item $h$ is bounded on the stabilizer $H(\alpha)$ of every essential
(non-degenerate if $S$ is non-orientable, see Appendix)
  simple closed curve (or scc) $\alpha$
in $S$.
\end{enumerate}
\end{thm}

Since we did not state the theorem in exactly this form in \cite{bf}
we will outline the proof in Section \ref{outline}.
Also, in \cite{bf} we discussed only the case when $S$ is orientable.
We explain how to modify the argument when $S$ is non-orientable.

\begin{cor}\label{slim}
Suppose $MCG(S)$ is not virtually abelian.  Let $$Q_1,\cdots,Q_m$$ be a
finite collection of subgroups of $MCG(S)$ such that each $Q_i$ is
either cyclic or the stabilizer of an essential (non-degenerate)
scc. Then
$$Q_1Q_2\cdots Q_m\neq MCG(S)$$
Indeed, any finite union of sets of the above form is a proper subset
of $MCG(S)$.
\end{cor}

Here the set on the left consists of compositions $\phi_1\cdots\phi_m$
for $\phi_i\in Q_i$. This statement generalizes the fact that $MCG(S)$
is not boundedly generated (by cyclic subgroups) \cite{FLM}.

In this note we will extend Theorem \ref{old} so that (2) includes a
larger class of subgroups of $MCG(S)$. The additional subgroups will
be called {\it D-subgroups} and they consist of mapping classes that
descend down a particular proper covering map $p:S\to S'$. The
motivation for seeking such an extension comes from the first order
theory of free groups, which also demands a consideration of
non-orientable surfaces. In this theory one is interested in
understanding the set of conjugacy classes of homomorphisms $f:H\to \F$
that extend to $\tilde f:G\to \F$, where $\F$ is a fixed nonabelian
free group, and $H\subset G$ are fixed finitely generated groups. For
example, consider the special case when $H$ is the fundamental group
of a closed surface $S$.  To glean the structure of this set it is
useful to consider its intersections with ``orbits'': Let $f:H\to\F$
be a fixed homomorphism and consider the set
$$E_f=\{\phi\in MCG(S)|f\phi \mbox{ extends to }G\}$$
With the help of Makanin-Razborov diagrams \cite{zlil} (see also
\cite{km1} and \cite{km2}) one shows that either $E_f=MCG(S)$ or $E_f$
is contained in a finite union of sets of the form
$$Q_1Q_2\cdots Q_m$$
where each $Q_i$ is either cyclic, or it is contained in the
stabilizer of a simple closed curve, or it is a D-subgroup. The
details of this are in \cite{bf:tarski}, where a sharper statement is
proved. Thus Corollary \ref{slim} implies (using our extension of
Theorem \ref{old}) that $E_f$ is either all of $MCG(S)$ or else it is
a ``slim'' subset of $MCG(S)$. 

This reasoning applies to groups $H$ other than surface groups by
considering abelian JSJ decompositions \cite{ds}, \cite{fp}. The most
interesting pieces of the decomposition are surfaces (with boundary)
and we are led to considering non-closed surfaces as well. It is
convenient to collapse boundary components to punctures. This does no
harm, since it only means that the quasi-homomorphisms we construct
will be 0 on the subgroup generated by Dehn twists in the boundary
components.

The crucial tool in the study of mapping class groups is that of the
{\it curve complex} associated to the surface. This complex was
originally defined by Harvey \cite{harvey}, it was used by Harer
\cite{harer} in his study of the homology of mapping class groups, and
its geometric aspects have been studied more recently thanks to the
celebrated theorem of Masur and Minsky which states that the curve
complex is $\delta$-hyperbolic. We review the definition and basic
properties of the curve complex in the Appendix, where we also outline
proofs of basic facts for the case of non-orientable surfaces.

\vskip 0.5cm
\noindent
{\bf Acknowledgement.}
The second author would like to thank Mustafa Korkmaz for useful
discussions.

\section{D-subgroups of mapping class groups}

Let $p:S\to S'$ be a covering map between two compact connected
surfaces. We will also assume that $p$ is not a homeomorphism. There
is a finite index subgroup $MCG_p(S')<MCG(S')$ consisting of (isotopy
classes of)
homeomorphisms $\phi':S'\to S'$ that admit a lift, i.e. a
homeomorphism $\phi:S\to S$ such that $p\phi=\phi' p$. When a lift
exists it may not be unique. Consider the group 
$$D(p)=\{\phi:S\to S|p\phi=\phi' p \mbox{ for some }\phi':S'\to S'\}/\mbox{isotopy}<MCG(S)$$
consisting of all possible lifts of elements of
$MCG_p(S')$. Equivalently, this is the group of mapping classes in
$MCG(S)$ that {\it descend} to $S'$ via $p$.

The group $D(p)$ is {\it commensurable} with $MCG_p(S')$ in the sense
that the group $$\tilde
D(p)=\{(\phi,\phi')|p\phi=\phi'p\}/\mbox{isotopy}$$ surjects to both
with finite kernel.

Any subgroup of $MCG(S)$ of the form $D(p)$ will be called a {\it
  D-subgroup} of $MCG(S)$. Such subgroups arise naturally in the study
  of the first order theory of free groups and they play a role
  similar to the stabilizer of a simple closed curve in $S$. 

\begin{lemma}\label{finite}
If $\chi(S)<0$ then there are only finitely many conjugacy classes of
$D$-subgroups of $MCG(S)$.
\end{lemma}

\begin{proof}
There are finitely many possible topological types of surfaces $S'$
that can appear in a covering map $p:S\to S'$.  Fix one such $S'$.
Unless $S'$ is a disk with 2 holes or a projective plane with 2 holes,
choose two tight simple closed curves $\alpha$ and $\beta$ that fill
$S'$ (this means that each complementary component is a disk or an
annulus that contains a boundary component of $S'$). If $S'$ is a disk
with 2 holes choose an immersed curve $\alpha$ with one
self-intersection point so that all 3 complementary components are
annuli containing one boundary component of $S'$. Similarly, when $S$
is a projective plane with 2 holes choose $\alpha$ so that it has two
self-intersection points and the complementary components are
annuli. We can view this as a kind of a cell structure on $S'$ except
that some of the cells have a hole in the interior.

If $p:S\to S'$ is a covering map we can pull back this ``cell
structure'' to $S$. The cells are labeled by the corresponding cells
in $S'$. There are only finitely many labeled cell complexes that
arise in this fashion. Now suppose that $p_1,p_2:S\to S'$ are two
covering maps with the associated labeled cell complexes isomorphic
via a homeomorphism that preserves labels. This yields $h:S\to S$ such
that $p_1=p_2h$, and this means that $D(p_1)$ and $D(p_2)$ are
conjugate in $MCG(S)$: $D(p_1)=h^{-1}D(p_2)h$.
\end{proof}

\section{Outline of proof of Theorem \ref{old}}\label{outline}

Suppose $S$ is orientable.
Since $MCG(S)$ is not virtually abelian, $S$ is not 
$S^2$ minus at most three points.
If $S$ is $S^2$ minus four points or $T^2$ minus at most one point, 
then $MCG(S)$ is virtually free. In those cases, the curve complex $X$
(see Appendix)
is not connected, and we need to argue 
differently. One way is to modify the definition of $X$ (see Apendix), 
then we can apply the same argument as follows. We omit the details.

The mapping class group $MCG(S)$ acts on the curve complex $X=X(S)$ which is
  hyperbolic
by the
  celebrated theorem of Masur-Minsky \cite{MM}. For finite oriented
  paths $w$, $\alpha$ in $X$ write $|\alpha|_w$ for the maximal number
  of non-overlapping translates of $w$ in $\alpha$, and by $|w|$
  denote the length of $w$. If $W$ is an integer with $0<W<|w|$, for
  any two vertices $x,y\in X$ define
$$c_{w,W}(x,y)=d(x,y)-\inf_\alpha (|\alpha|-W|\alpha|_w)$$
with $\alpha$ ranging over all paths from $x$ to $y$. Finally, define
$$h_w:MCG(S)\to\R$$
by
$$h_w(g)=c_{w,W}(x_0,g(x_0))-c_{w^{-1},W}(x_0,g(x_0))$$
where $x_0\in X$ is a fixed base vertex.

Then $h_w$ is a quasi-homomorphism \cite[Proposition 3.10]{koji}. By
construction, $h_w$ is bounded on the stabilizer of every vertex $x$ of
$X$ (i.e. the stabilizer of a scc in $MCG(S)$). Indeed, $0\leq
c_{w,W}(x_0,g(x_0))\leq d(x_0,g(x_0))\leq 2d(x_0,x)$, so $|h_w(g)|\leq
2d(x_0,x)$ for any $g\in MCG(S)$ that fixes $x$. 

In \cite[Proposition 11]{bf} we showed that the action of $MCG(S)$ on
$X$ satisfies a certain technical condition called WPD (see Appendix). In the
presence of this condition, the following two statements about
hyperbolic elements $g_1,g_2$ are
equivalent \cite[Proposition 6]{bf}:
\begin{itemize}
\item some positive powers $g_1^{n_1}$ and $g_2^{n_2}$ are conjugate,
  and
\item for any quasi-axes $\ell_1$ and $\ell_2$ of $g_1$ and $g_2$
  respectively there is a constant $C>0$ such that for any $L>0$ there
  are segments $\Sigma_i\subset\ell_i$ of length $L$ so that a
  translate of $\Sigma_1$ is contained in the $C$-neighborhood of
  $\Sigma_2$ and the orientations of the two segments are
  parallel.
\end{itemize}
Write $g_1\sim g_2$ if the two statements hold. We then constructed in
\cite[Proposition 2]{bf} an infinite sequence $f_1,f_2,\cdots\in
MCG(S)$ of hyperbolic elements such that
\begin{itemize}
\item $f_i\not\sim f_i^{-1}$ for all $i$, and
\item $f_i\not\sim f_j^{\pm 1}$ for $i\neq j$.
\end{itemize}
If $C=<\phi>$ is a cyclic subgroup of $MCG(S)$ generated by a
  hyperbolic element $\phi$ then $\phi\sim f_i^{\pm
  1}$ for at most one $i$. So if we are given a finite collection of
  cyclic subgroups generated by hyperbolic elements we may assume
  (after removing a finite subset of the $f_i$'s) that
  they do not contain any conjugates of any nontrivial powers of the
  $f_i$'s. Now, as in the proof of \cite[Theorem 1]{bf}, inductively
  choose paths $w_1,w_2,\cdots$ where $w_i$ is a long segment along a
  quasi-axis of $f_i$ so that $h_{w_i,W}:MCG(S)\to\R$ is unbounded on
  $<f_i>$ but bounded on the given cyclic subgroups and on
  $<f_1>,\cdots,<f_{i-1}>$ (when the cyclic subgroup is generated by
  an element that is not hyperbolic, then any $h_{w,W}$ is bounded on
  it).
Then the quasi-homomorphisms $h_{w_i,W}$ give an infinite linearly
  independent collection in $QH(MCG(S))$ and by construction they are
  all bounded on the subgroups in (1) and (2).

Suppose $S$ is non-orientable.
Since $MCG(S)$ is not virtually abelian,
$S$ is not ${\Bbb R}P^2$ minus at most two points
($MCG$ is finite, \cite{kork}), 
nor Klein bottle minus at most one point
($MCG$ is finite (no puncture) or virtually ${\Bbb Z}$ (one puncture)
\cite{stu}).
As we show in the appendix, the curve complex $X$ is connected, and 
delta-hyperbolic, and the action of $MCG$ is WPD. Therefore the same 
argument applies.

\section{Extension of Theorem \ref{old}}

\begin{lemma}\label{quasi-convex}
Let $p:S\to S'$ be a covering map between compact connected surfaces
of negative Euler characteristic, possibly with boundary, possibly
nonorientable. Then
the image of
$$p^*:X(S')\to X(S)$$ is a quasi-convex subset of $X(S)$.
Namely, there exists a constant $P$ such that 
a geodesic (in $X(S)$) between any two points in $p^{*}(X(S'))$
is in the $P$-neighborhood of $p^{*}(X(S'))$.
\end{lemma}

This lemma complements the result of Rafi and Schleimer \cite{RS} that
$p^*$ is a quasi-isometric embedding.

\begin{proof}
We use terminology and results from Appendix.
Let $a,b$ be two vertices of $X(S')$ and let $c$ be a scc as in Lemma
\ref{1}. Likewise, let $\tilde a,\tilde b$ be the preimages of $a,b$
in $S$, viewed as multi-curves. Let $\tilde c$ be the preimage of
$c$. Since the length $l_{ab}(c)$ of $c$ is bounded above by
$R\sqrt{I(a,b)}$ it follows that $$l_{\tilde a\tilde b}(\tilde c)\leq
|p|R\sqrt{I(a,b)}=\sqrt{|p|}R\sqrt{I(\tilde a,\tilde b)}$$
and similarly
$$m_{\tilde a\tilde b}(\tilde c)=\sup\frac{I(\tilde c,\tilde
  x)}{l_{\tilde a\tilde b}(\tilde x)}=\sup
  \frac{I(c,x)}{l_{ab}(x)}=m'_{ab}(c)\le
  R/\sqrt{I(a,b)}=\sqrt{|p|}R/\sqrt{I(\tilde a,\tilde b)}$$ where $\tilde x$
  run over scc's in $S$ and $x$ denotes $p(\tilde x)$ (the modulus can
  be computed by running the sup over all immersed curves, see Remark
  \ref{modulus}). Therefore, $\tilde c$ is within a bounded distance
  from $Mid(\tilde a,\tilde b)$ (see Remark \ref{m2l}). A similar
  argument applies to $Mid(\tilde a,\tilde b;\frac pq)$.
\end{proof}

In the orientable case the following proposition is a consequence of
\cite[Theorem 2]{MR1214233}. 

\begin{prop}\label{generic pA}
Let $S$ be a compact surface that admits pseudo-Anosov homeomorphisms.
Then there is a pseudo-Anosov homeomorphism $\phi:S\to S$ such that no
nontrivial power of $\phi$ (is conjugate within $MCG_{\pm}(S)$ to a
homeomorphism that) descends down a proper covering $S\to S'$.
\end{prop}

\begin{proof}
In this proof it is convenient to collapse boundary components to
punctures. If $S$ is a 4 times punctured sphere, a once punctured
torus, or a once punctured $\#_1^3\R P^2$ then $S$ admits a
pseudo-Anosov homeomorphism, and such a homeomorphism cannot descend
down a proper covering map since simpler surfaces do not admit such
homeomorphisms. (For the first two cases $\chi=-1$ so the
assertion is obvious. For the last case this takes a thought since
$\chi=-2$, so conceivably there is a double cover. But the covered
surface would have to have one puncture and under a double cover a
single puncture must lift to two punctures for homological reasons.)

In all other cases we will construct $\phi:S\to S$ such that the
stable and unstable foliations of $\phi$ have a unique singularity with $>2$
prongs. This will imply that $\phi$ and its powers do not descend down
a proper covering map. Note that if we have such a map $\phi$ for a surface
$S$ with $p$ punctures, then we also have it for the surface $S$ with
$>p$ punctures, by taking an appropriate power of $\phi$ and declaring
fixed points to be punctures. Thus it suffices to construct $\phi$ in
the ``minimal'' cases.

Recall Thurston's construction of pseudoAnosov homeomorphisms of a
surface $S$ \cite{thurston}. Let $\alpha$ and $\beta$ be two tight
collections of pairwise disjoint essential 2-sided simple closed
curves in $S$ such that $\alpha\cup\beta$ fills $S$. By composing Dehn
twists in these curves one obtains a pseudo-Anosov homeomorphism
$\phi:S\to S$ whose (un)stable foliation has one singularity for every
complementary component of $\alpha\cup\beta$. If the component is a
$2k$-gon, then the singularity is $k$-pronged (so if $k=2$ it's not a
singularity, and if $k=1$ there must be a puncture in this component).
The proof is now contained in Figures 1-8.
%\begin{minipage}{5cm}
\floatstyle{ruled}
\newfloat{type}{b}{cap}

\begin{figure}%[h!]
\centerline{\includegraphics[scale=0.5]{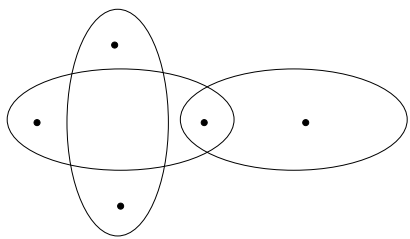}}
\caption{$S^2$ with 5 punctures. Each puncture is a 1-pronged
  singularity, and there is a 3-pronged singularity.}
\end{figure}
%\end{minipage}
%\begin{minipage}{5cm}
\begin{figure}%[h!]
\centerline{\includegraphics[scale=0.5]{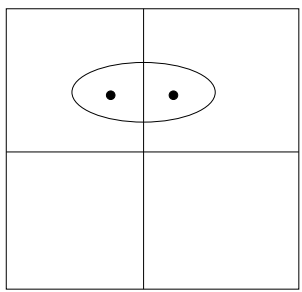}}
\caption{$T^2$ with 2 punctures, pictured as the square with opposite
  sides identified. The curves are the meridian and the longitude as
  well as a curve surrounding the two punctures. Each puncture is a 1-pronged
  singularity, and there is a 4-pronged singularity.}
\end{figure}
%\end{minipage}
\begin{figure}
\centerline{\includegraphics[scale=0.4]{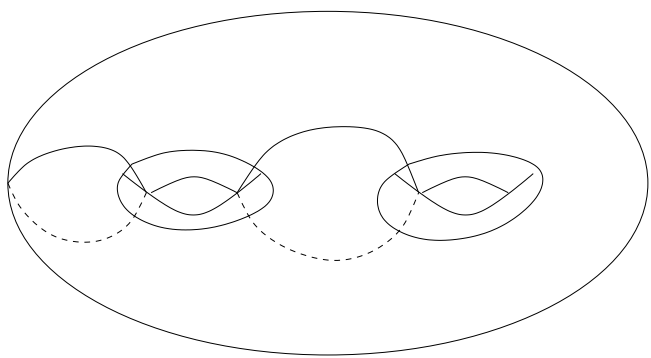}}
\caption{Genus 2 surface with a single 6-pronged singularity. In a
  similar way we get a genus $g\geq 2$ surface with a single
  $4g-2$-pronged singularity.}
\end{figure}
\begin{figure}
\centerline{\includegraphics[scale=0.4]{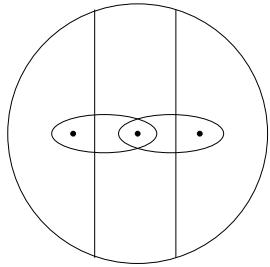}}
\caption{$\R P^2$ with 3 punctures, pictured as a disk with antipodal
  points on the boundary identified. Each puncture is a 1-pronged
  singularity, and there is a 3-pronged singularity.}
\end{figure}
\begin{figure}
\centerline{\includegraphics[scale=0.4]{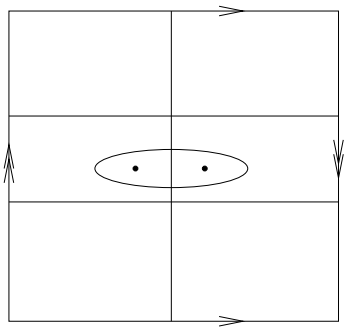}}
\caption{The Klein bottle $K$ with 2 punctures. Each puncture is a 1-pronged
  singularity, and there is a 4-pronged singularity.}
\end{figure}
\begin{figure}
\centerline{\includegraphics[scale=0.4]{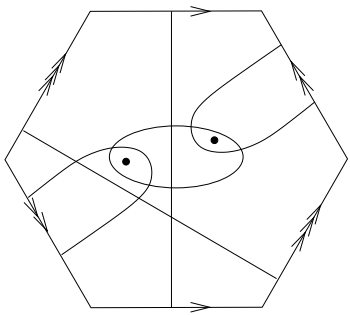}}
\caption{Surface with $\chi=-1$ and 2 punctures. Each puncture is a 1-pronged
  singularity, and there is a 6-pronged singularity.}
\end{figure}
\begin{figure}
\centerline{\includegraphics[scale=0.4]{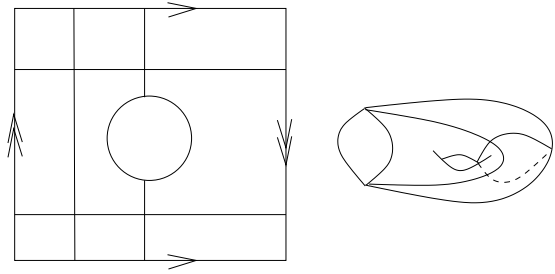}}
\caption{$K\#T^2$ with a single 4-pronged singularity. In a similar way
  we get $K\#\#_1^nT^2$ with a single $4n$-pronged singularity.}
\end{figure}
\begin{figure}
\centerline{\includegraphics[scale=0.5]{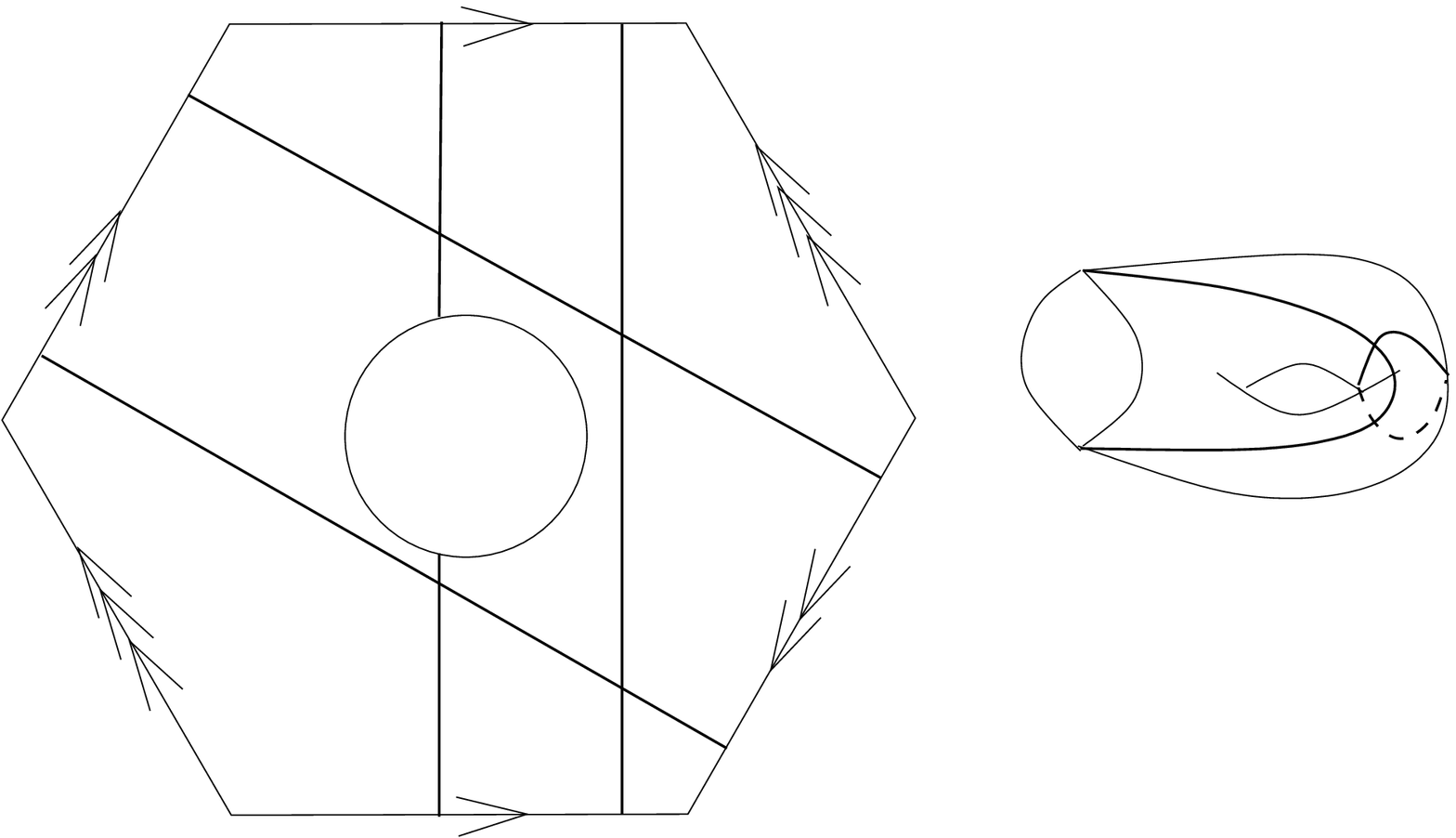}}
\caption{$\#_1^3\R P^2\#T^2$ with a single 6-pronged singularity. In a
  similar way we get $\#_1^3\R P^2\#\#_1^nT^2$ with a single
  $4n+2$-pronged singularity.}
\end{figure}
\end{proof}

\begin{cor}\label{one map}
Let $\phi:S\to S$ be as in Proposition \ref{generic pA} and let $g$ belong
to the
centralizer $Z(\phi)$ of $<\phi>$ in $MCG(S)$. If $g$ has finite order
and is realized as an isometry of a hyperbolic structure on $S$ then
this isometry necessarily has a fixed point.
\end{cor}

\begin{proof}
Fix a hyperbolic structure on $S$ such that $g$ is realized by an
isometry $S\to S$. Let $\Lambda$ be the stable geodesic lamination
associated with $\phi$. Then $g(\Lambda)$ is the stable geodesic
lamination associated with $g\phi g^{-1}=\phi$ and so
$g(\Lambda)=\Lambda$. By construction, there is $k>2$ and a unique
complementary component of $\Lambda$ which is a $k$-gon. Therefore $g$
maps this component to itself and has a fixed point.
\end{proof}

We will need a small generalization of this.

\begin{cor}\label{many maps}
Suppose $\phi:S\to S$ is a pseudo-Anosov homeomorphism whose stable
foliation has precisely $n$ $k$-pronged singularities for a certain
$k>2$. Let $\Delta$ be a finite subgroup of the centralizer $Z(\phi)$,
and suppose that $\Delta$ is realized as a group of isometries of a
hyperbolic metric on $S$. If $|\Delta|>n$ then $\Delta$ does not act
freely on $S$.
\end{cor}

\begin{proof}
Now $\Delta$ permutes the complementary $k$-gons and $|\Delta|>n$
ensures that the stabilizer of a $k$-gon is nontrivial, and fixes a
point.
\end{proof}

\begin{cor}\label{short overlaps}
Fix a homeomorphism $\phi:S\to S$ as in Proposition \ref{generic
pA}. Let $\ell$ be a quasi-axis of $\phi$ in the curve complex $X(S)$. Then
for every proper covering map $p:S\to S'$ and every $B>0$ there is $L>0$ such
that no translate of $\ell$ has a segment of length $L$ that is
contained in the $B$-neighborhood of $p^*(X(S'))$.
\end{cor}

\begin{proof}
First suppose the covering map $p':S\to S'$ is regular. Let $\Delta$
be the deck group. Then $\Delta$ acts simplicially on the curve
complex $X(S)$ and the fixed set is precisely $p^*(X(S'))$ (it is a
subcomplex of the barycentric subdivision of $X(S)$). Now suppose that
some translate $h(\ell)$ has a segment of length $L$ contained in the
$B$-neighborhood of $p^*(X(S'))$. Then for every $g\in\Delta$ we have
that $gh(\ell)$ and $h(\ell)$ have segments of length $L$ that are
within $2B$ of each other. From \cite{bf} we conclude that when $L$ is
sufficiently large there is $N>0$ such that $gh\phi^N
h^{-1}g^{-1}=h\phi^N h^{-1}$. This means that $g$ is in the
centralizer of $h\phi^N h^{-1}$. But then Corollary \ref{one map} (applied to
$h\phi^N h^{-1}$ and $g$) implies that when $g$ acts as an
isometry of a hyperbolic structure on $S$ then it fixes a point. But by
construction, $g$ acts without fixed points.

Now consider the general case, when $p:S\to S'$ is not regular. Let
$q:\tilde S\to S$ be a finite cover with $\tilde S$ connected so that
$pq$ (and hence $q$) is regular. Let $\tilde\Delta$ be the deck group
of $pq$ and $\Delta\subset\tilde\Delta$ the deck group of $q$. Again
suppose that $h(\ell)$ contains a segment of length $L$ that is inside
the $B$-neighborhood of $p^*(X(S'))$, for some $h\in MCG(S)$. For some
$M>0$ the homeomorphism $h\phi^M h^{-1}$ lifts to $\tilde S$. Fix one
such lift $\tilde \tau:\tilde S\to\tilde S$ and note that 
\begin{itemize}
\item the stable foliation of $\tau$ has precisely $|q|$ ($=$ number of sheets)
  singularities with $k$ prongs, for a certain $k>2$, and
\item $q^*(\ell)\subset X(\tilde S)$ is a quasi-axis of $\tau$.
\end{itemize}
Of course, a segment of length $L$ in $q^*(\ell)\subset
X(\tilde S)$ is also contained in the $B$-neighborhood of
$q^*(p^*(X(S')))$. It now follows exactly as above that for some $N>0$
the deck group $\tilde\Delta$ is contained in the centralizer of
$\tilde\tau^N$. Corollary \ref{many maps} then says that $\tilde\Delta$
cannot act freely as an isometry group on $\tilde S$, contradiction.
\end{proof}

\begin{thm}\label{new}
Let $S$ be a compact surface such that $MCG(S)$ is not virtually
abelian. Then $QH(MCG(S))$ is infinite-dimensional. Moreover, for
every finite collection of cyclic subgroups $C_1,\cdots,C_k\subset
MCG(S)$ there is an infinite-dimensional subspace of $QH(MCG(S))$ such
that each representative quasi-homomorphism $h$ of any element in this
subspace satisfies:
\begin{enumerate}[(1)]
\item $h$ is bounded on each $C_i$, 
\item $h$ is bounded on the stabilizer $H(\alpha)$ of every essential
(non-degenerate when $S$ is non-orientable)
  simple closed curve $\alpha$
in $S$, and
\item $h$ is bounded on every D-subgroup $D(p)$.
\end{enumerate}
\end{thm}

\begin{proof}
The strategy of proof is the same as that of Theorem \ref{old} but
certain additional things need to be arranged.

Let $\phi:S\to S$ be the
pseudo-Anosov homeomorphism from Proposition \ref{generic pA}. Also
choose another pseudo-Anosov homeomorphism $\psi:S\to S$ so that
$\phi,\psi$ do not have isotopic nontrivial powers. The subgroup
$F\subset MCG(S)$ generated by high powers $\phi^N,\psi^N$ of $\phi$
and $\psi$ is a nonabelian free group and moreover an $F$-equivariant
map $j:T\to X$ is a quasi-isometric embedding, where $T$ is the Cayle
tree of $F$ with respect to $\phi^N,\psi^N$. 
Thus for every $1\neq f\in F$ we obtain a quasi-axis for
$f$ by taking the axis of $f$ in $T$ and mapping it to $X$ by
$j$. There are constants $K_0,C_0$ so that in this way we obtain a
$(K_0,C_0)$-quasi-axis regardless of which $1\neq f\in F$ we
chose. All paths $w$ used in the construction of a quasi-homomorphism
$h_{w,W}$ will be $j$-images of edge-paths in $T$, and in particular
they are all $(K_1,C_1)$-quasi-geodesic for a fixed $K_1,C_1$. It now
follows that every minimizer path is a $(K,C)$-quasi-geodesic for a
fixed $K,C$ (see the proof of Theorem
\ref{old} given in Section \ref{outline}).

Let $D_1,\cdots,D_m$ be representatives of conjugacy classes of
$D$-subgroups of $MCG(S)$ (see Lemma \ref{finite}). If a
quasi-homomorphism is bounded on each $D_k$ then it is bounded on
every $D$-subgroup. 
Say $D_k=D(p_k:S\to S_k')$ and let $A_k=p_k^*(X(S_k'))\subset
X(S)$. Then $A_k$ is quasi-convex in $X(S)$ (see Lemma
\ref{quasi-convex}). Choose $B>0$ such that any
$(K,C)$-quasi-geodesic joining two points of $A_k$ is contained in the
$B$-neighborhood of $A_k$, $k=1,\cdots,m$ and choose $L>0$ such that
no segment in $\ell$ (the quasi-axis of $\phi$) can be translated into
the $B$-neighborhood of any $A_k$ (see Corollary \ref{short
  overlaps}). 

Now choose $f_1,f_2,\cdots\in F$ so that in addition to $f_i\not\sim
f_j^{\pm 1}$ for $i\neq j$ and $f_i\not\sim f_j^{-1}$, the quasi-axis
of each $f_i$ contains a path in a translate of $\ell$ of length
$>L$. Finally, choose a path $w_i$ in the quasi-axis of $f_i$ so that
in addition to the previous requirements $w_i$ contains a segment of
length $>L$ which is a translate of a segment in $\ell$. It follows
that no $w_i$ can be translated into the $B$-neighborhood of any
$A_k$. In particular, $(K,C)$-quasi-geodesics joining points of $A_k$
do not contain copies of $w_i$. If the base vertex $x_0$ belongs to
$A_k$ then we deduce that $h_{w_i,W}(g)=0$ for every $g\in MCG(S)$
that leaves $A_k$ invariant. When $x_0$ is moved, $h_{w,W}$ changes by
a bounded amount.
\end{proof}

\begin{cor}\label{maincor}
Suppose $MCG(S)$ is not virtually abelian.
Let $Q_1,\cdots,Q_m$ be a finite collection of subgroups of $MCG(S)$
such that each $Q_i$ is either cyclic, or the stabilizer of a 
non-degenerate scc, or
a $D$-subgroup. Then
$$Q_1Q_2\cdots Q_m\neq MCG(S)$$
Indeed, any finite union of sets of the above form is a proper subset
of $MCG(S)$.
\end{cor}

\section{Appendix}

\renewcommand{\Sigma}{S}

In this Appendix we discuss the curve complex $X$ of a compact 
surface $S$, which is orientable or non-orientable.
Masur and Minsky \cite{MM} have shown that $X$ is Gromov-hyperbolic
when $S$ is orientable.
Bowditch \cite{bowditch:hyp} has given another argument for that fact
using CAT(0) geometry.
We follow his argument and explain that it also
applies to non-orientable surfaces with minor changes.

The mapping class group $MCG(S)$ acts on $X$ by isometries.
When $\Sigma$ is orientable, 
Masur and Minsky  have shown that an element in $MCG(\Sigma)$
acts by a hyperbolic isometry on $X$ if and only if 
it is pseudo-Anosov. It is shown that the action 
is weakly properly discontinuous \cite{bf}. 
We explain that those facts remain true when $\Sigma$ is non-orientable.
The argument does not require any change.

\subsection{The curve complex}

Let $\Sigma$ be a closed surface equipped with $n$ punctures
(i.e. distinguished points).  

When $\Sigma$ is orientable, let $X=X(\Sigma,n)$ be the simplicial
complex whose vertices are isotopy classes of {\it essential} scc's
(i.e. those that don't bound a disk or a
punctured disk). 

When $\Sigma$ is non-orientable, vertices
are isotopy classes of scc's {\it non-degenerate}
(does not bound a disk, a punctured disk,
or a M\"obius band).
If $c$ bounds a M\"obius band, the core curve
of the M\"obius band is non-degenerate.

In either case, a collection of vertices of $X$ bounds a simplex
provided representative curves can be found that are pairwise
disjoint.

Denote by $I(a,b)$ the minimal intersection number between
representatives. It is attained when there are no bigons in the
complement. 

When $\Sigma=S^2$ and $n\leq 3$ then $X=\emptyset$. When $\Sigma=S^2$
and $n=4$, or $\Sigma=T^2$ and $n\leq 1$, then $X$ is a discrete set
(but see Remark \ref{farey} below). So we will assume that $n\geq 5$
when $\Sigma=S^2$ and $n\geq 2$ when $\Sigma=T^2$.  Similarly, we will
assume $n\ge 3$ when $\Sigma$ is $\R P^2$, and we assume $n\ge 2$ when
$\Sigma$ is the Klein bottle.

\begin{lemma}
$X$ is connected.
\end{lemma}

\begin{proof}
Let $a,b$ be two vertices. We will connect them by an edge-path.

If $I(a,b)=0$ then there is an edge between them. Next, we consider
the case that both $a$ and $b$ are 2-sided.
If $I(a,b)=1$ then the regular neighborhood $N$ of $a\cup b$ is a torus
with one boundary component $c$. This is an essential curve since
otherwise we are in a torus with $n\leq 1$. If $c$ bounds a M\"obius
band, replace $c$ by the core curve. Thus $c$ is connected to
both $a$ and $b$.

As the next case, suppose we can find two consecutive intersection
points on $a$ that have the same sign (orient $a$ and transversely
orient $b$, so signs make sense). Form a curve $c$ as the union
of the arc in $a$ that joins these two points and has no intersection
points in the interior, and of an arc in $b$ joining these two
points. Then either $I(b,c)=1$ or $I(b,c)=0$ and $c$ is one-sided. In
either case, $c$ is non-degenerate. Since $I(a,c)<I(a,b)$ the claim
now follows by induction.

Now suppose that we have two intersection points of the same sign on
$a$ with one other intersection point in between. Then we can choose
an arc on $b$ connecting the two points that does not intersect the
interior of the arc on $a$, so the same argument works.

Finally, if $a$ and $b$ intersect in two points of opposite signs,
then the regular neighborhood $N$ is a sphere with 4 boundary
components. At least one boundary component must be essential, and if
it is degenerate replace it by the core of the M\"obius band it bounds, so
$dist(a,b)=2$.

Now suppose at least one of $a,b$, say $b$, is 1-sided. If $I(a,b)=1$
then the regular neighborhood $N$ of $a\cup b$ is either $\R P^2$ with
two boundary components or the Klein bottle with one boundary
component. In either case, at least one boundary component $c$ must be
essential. If $c$ is degenerate, replace it by the core of the
M\"obius band it bounds. Thus $d(a,b)=2$.

Now suppose $I(a,b)>1$.
Fix two consecutive intersection points on $a$.  Using the
arc in $a$ between the two points (and contains no other intersection
points) and an arc in $b$ between the two points (there are two such
arcs in $b$), form a scc $b_1$. Using the other arc in $b$, form a scc
$b_2$.  By construction, $I(b_1,b), I(b_2,b) \le 1$, and also
$I(b_1,a), I(b_2,a) < I(b,a)$.

Note that $b_1+b_2=b$ in $H^1(\Sigma,{\Bbb Z}_2)$.
A scc is one-sided if and only if its square
is non-trivial in $H^2(\Sigma,{\Bbb Z}_2)$.
By our assumption, $b^2$ is non-trivial 
in $H^2(\Sigma,{\Bbb Z}_2)$.
Therefore, one of $b_1,b_2$, say $b_1$, is 
one-sided. It follows that $b_1$ is 
non-degenerate, which gives a vertex
in $X$, so we are done by induction. 
\end{proof}

Higher connectivity of the complex of curves for orientable 
and non-orientable $S$ is 
discussed in \cite{ivanov} (for example, Th 2.6).

\begin{remark}
Induction shows easily that $dist(a,b)\leq 2I(a,b)$. This argument is
essentially contained in \cite[Lemma 2]{Lickorish}.
Actually, the argument can be improved.
In the generic case above we can choose one of two possible arcs in
$b$ connecting the two intersection points. We choose the one with
fewer intersection points. This guarantees that $I(a,c)<I(a,b)-\big[
  \frac{I(a,b)-1}2\big]$. Induction then shows that $dist(a,b)\leq
O(\log I(a,b))$.
\end{remark}

\begin{remark}
There is no lower bound to $dist(a,b)$ in terms of $I(a,b)$. Just
think about two curves that intersect a lot but don't fill. In fact,
it is nontrivial (and involves Nielsen-Thurston theory) to show that
$X$ has infinite diameter.
\end{remark}

\begin{remark}\label{farey}
When $\Sigma=T^2$ and $n=0$ or 1, it is natural to change the
definition and connect $a$ to $b$ when $I(a,b)=1$. The resulting graph
is classically known as the Farey graph -- it is the 1-skeleton of the
standard tesselation of ${\mathbb H}^2$ by ideal triangles, including
the vertices at infinity. There is
classical number theory (continued fractions) associated with this
graph -- see the article by C. Series \cite{MR810563}. This graph is also
hyperbolic. Exercise: Give a direct proof of this fact.
\end{remark}

\subsection{Hyperbolicity of the curve complex}
Masur-Minsky \cite{MM} showed that the curve complex of an orientable
surface is delta-hyperbolic.  Bowditch \cite{bowditch:hyp} gave
another proof of this fact.  We outline his argument and explain that
his argument applies to non-orientable surfaces as well.

By a {\it curve system} we mean a collection of essential, pairwise
disjoint scc's (we allow parallel curves).

We will think of $\Sigma$ as a compact surface with punctures.
Let $a,b$ be two curve systems on the surface $\Sigma$. Realize them
so that they intersect minimally (every bigon contains a puncture). We
will assume that $a\cup b$ is {\it filling}, i.e. each disk in the
complement contains a puncture, and we will construct a particular
Euclidean metric on $\Sigma$ with cone singularities. All cone angles
will be multiples of $\pi$ and if an angle equals $\pi$ then the point
is one of the punctures. Realize $a$ and $b$ in $\Sigma$ so that they
intersect minimally (every bigon contains a puncture).
Take a collection of squares, one for each
intersection point, and glue them in the obvious way so that $a\cup b$
is the dual 1-skeleton of the resulting 2-complex. This complex is
$\Sigma$. Also, the horizontal and vertical foliations on the squares
match up and give a pair of (measured) foliations on $\Sigma$. The
area of $\Sigma$ equals the intersection number $I(a,b)$. By $l_S(c)$
denote the length (with respect to this metric) of a shortest curve
representing $c$.

Now set:
\begin{itemize}
\item the {\it length} of a curve $c$ is
$l(c)=l_{ab}(c)=I(a,c)+I(b,c)$,
\item the {\it modulus} of $c$ is
  $m(c)=m_{ab}(c)=\sup\{\frac{I(c,x)}{l(x)}|x\mbox{ is a curve}\}$.
\end{itemize}

We will restrict ourselves to the case when all curves in $a$ and $b$
are 2-sided. The reason is that we want to be able to talk e.g. about the
curve system $Na$ for $N=1,2,3,\cdots$, and this can be taken to mean
``take $N$ parallel copies of $a$''. If $a$ were 1-sided, we could
only do this formally.

\begin{lemma}
Suppose each curve in $a$ and $b$ is 2-sided. Then 
$$\frac 1{\sqrt 2}l_S(c)\le l_{ab}(c)\le 2l_S(c)$$
\end{lemma}

\begin{proof}
First suppose that the geodesic representing $c$ is transverse to $a$
and $b$. The first inequality follows from the fact that each
complementary component of $a\cup b$ has diameter $\sqrt 2$, and the
second from the fact that the distance between two components of $a$
(and of $b$)
is $\ge 1$, so that $I(a,c)\le l_S(c)$ and $I(b,c)\le l_S(c)$.

Now suppose that $c$ is one of the components of $a$ (say). Then
$l_S(c)=I(b,c)=I(a,c)+I(b,c)=l_{ab}(c)$. 
\end{proof}

\begin{remark}
If $a$ is 1-sided curve and $c=a$ then
$l_S(c)=I(b,c)$ while $l_{ab}(c)=I(a,c)+I(b,c)=1+I(b,c)$, so the lemma
would still hold; however, it would fail if we allow formal scaling
for $c=Na$.
\end{remark}

\begin{lemma}\label{1}
Assume that if $\Sigma$ is $S^2$ then $n\ge 4$,
if $\Sigma=\R P^2$ then $n\ge 3$ and if $\Sigma$ is the Klein bottle
then $n\geq 2$.
There is a constant $R=R(\Sigma,n)$ such that for any filling pair $a,b$
of 2-sided curve systems
there is a curve $c$ with length bounded above by
$R\sqrt{I(a,b)}$ and modulus bounded above by $R/\sqrt{I(a,b)}$. The
curve $c$ is 2-sided, essential, possibly degenerate.
\end{lemma}

\begin{proof}(Outline)
We follow Bowditch \cite{bowditch:hyp}.
The argument requires only minor changes when $\Sigma$ is 
non-orientable. 

{\bf Step 1.}
We will find an essential
scc $c$ with length $l_S(c)$
bounded by $C\sqrt{Area(\Sigma)}$ and with an annular
(possibly Moebius band when $c$ is one-sided) neighborhood
around it of width at least $\sqrt{Area(\Sigma)}/C$ on each side
($C$ depends only on $\Sigma, n$).
This will
imply that for each $x$ we have $l_\Sigma(x)\geq
I(x,c)  \sqrt{Area(\Sigma)}/C$ and this implies that $m_{ab}(c)\leq
R/\sqrt{I(a,b)}$. Thus if $c$ is 2-sided we are done, and if $c$ is
1-sided replace it by the boundary of a regular neighborhood.

{\bf Step 2.} Quadratic isoperimetric inequality holds in
$\Sigma$. This means that if $D$ is a disk in $\Sigma$ that contains
at most one puncture, then $Area(D)\leq C\ length(\partial D)^2$ where
$C$ can be taken to be twice the constant for euclidean plane. Indeed,
when $D$ has no punctures this follows from the CAT(0) theory, and
when there is one puncture, pass to the double cover. Note that the
isoperimetric inequality fails when there are two punctures in $D$
(think of the double of an infinite half-strip). 

{\bf Step 3.} A {\it spine} in $\Sigma$ is a graph
$\Gamma\subset\Sigma\setminus\{\mbox{ punctures}\}$ such that
inclusion is $\pi_1$-surjective. We claim that there is a number
$\eta_1>0$ (as a function of the number of punctures) such that the
length of any spine satisfies the isoperimetric inequality (II)
$$length(\Gamma)\geq \eta_1\ \sqrt{Area(\Sigma)}.$$ It is enough to
prove this for ``minimal spines'', i.e. those that are no longer
spines after an edge is removed. When $\Sigma$ has no punctures, a
minimal spine has one complementary component, which is a disk, and
the standard II says:
$$Area(\Sigma)\leq C\ length(\partial D)^2=\frac
C4 length(\Gamma)^2$$ which proves the claim. When there are $n$
punctures, the complement consists of $n$ disks $D_i$ containing one
puncture each. From II we see that
$$Area(\Sigma)=\sum Area(D_i)\leq C\sum length(\partial D_i)^2$$
and so at least one $D_i$ must have long boundary, implying that
$\Gamma$ is long.

{\bf Step 4.} Let $N$ be the length of a longest chain of connected
incompressible, pairwise non-isotopic subsurfaces of $\Sigma\setminus
\{ \mbox{ punctures }\}$ and set $\eta_0=\eta_1/(100+2N)$. We claim
that there are two non-degenerate simple closed curves $x,y$ in
$\Sigma\setminus \{ \mbox{ punctures }\}$ whose  distance is
$\geq \eta_0\sqrt{Area(\Sigma)}$. 

For a curve $a$ denote by $length'(a)$ the length of $a$ when $a$ is
2-sided, and $\frac 32 length(a)$ when $a$ is 1-sided.  Let $x$ be a
non-degenerate scc in $\Sigma\setminus \{ \mbox{ punctures }\}$ with
$length'(x)$ within $\eta_0$ of $\inf\{length'(a)|a\mbox{ is
non-degenerate}\}$. We will show that there is a non-degenerate scc
$y$ outside the $\eta_0\sqrt{Area(\Sigma)}$-neighborhood of
$x$. Suppose not and consider the growing family of neighborhoods
$N_t$, $0\leq t\leq \eta_0\sqrt{Area(\Sigma)}$ of $x$. At $\leq N$
times the topology of the minimal incompressible subsurface that
contains $N_t$ changes. Each critical point can be connected to $x$ by
two arcs of length $\leq \eta_0\sqrt{Area(\Sigma)}$ so that the
obtained 1-complex $\Gamma$ consisting of these (disjoint, except for
endpoints) arcs and of $x$ carries $\pi_1$. The last surface is a
spine by our assumption, so this 1-complex is a spine. Thus,
$$length(x)+2N\eta_0\sqrt{Area(\Sigma)}\geq
\eta_1\sqrt{Area(\Sigma)}$$ and we deduce $$length(x)\geq
100\eta_0\sqrt{Area(\Sigma)}$$ Now it's not hard to find a curve in
$\Gamma$ much shorter than $x$.
More explicitly, we have the following cases:
\begin{itemize}
\item $x$ is 1-sided. If $\Gamma$ is a circle, then $\Sigma$ is 
$\R P^2$ with $\leq 1$ puncture. If there is a short arc in $\Gamma$
whose endpoints are far apart along $x$, then this arc together with
one of the two arcs of $x$ forms a scc $x'$ with $I(x,x')=1$ which is
much shorter than $x$. In particular, $x'$ is non-degenerate and
1-sided. If $\Gamma$ is the union of $x$ and only one short arc, then
$\Sigma$ is either $\R P^2$ with $\leq 2$ punctures or the Klein
bottle with $\le 1$ puncture. Now suppose that there are $\ge 2$ short
arcs in $\Gamma$ and each has endpoints close to each other along
$x$. For each such short arc consider the scc obtained by adding the
short subarc of $x$. If one such scc is nondegenerate then we have a
contradiction to the choice of $x$. Otherwise each such scc bounds a
punctured disk or a M\"obius band. In the latter case, the core would
be a better choice for $x$. If two of these curves bound punctured
disks then their ``connect sum'' along an arc of $x$ produces a curve
$x'$ which bounds a twice punctured disk and whose length is at most
slightly bigger than that of $x$. Thus $length'(x')\le \frac 23
length'(x)+\epsilon$ and $x'$ is a better choice than $x$.
\item $x$ is 2-sided and nonseparating. 
Then there is a short arc in $\Gamma$ that joins $x$ to itself from
  opposite sides. Let $x'$ be the curve formed by this arc together
  with the shorter of the two subarcs of $x$. Then $x'$ is much
  shorter than $x$ (even $\frac 32 length(x')$ is much shorter than
  $x$ in case $x'$ is 1-sided) and is non-degenerate because $I(x,x')=1$.
\item Now suppose $x$ is separating.
In particular, $x$ is 2-sided. If there is a short arc $\gamma$ in
$\Gamma$ such that $x\cup\gamma$ has a nonorientable regular
neighborhood, then $\gamma$ union the shorter of the two subarcs in
$x$ contradicts the choice of $x$. 
Similarly, if there is $\gamma$
whose endpoints are far apart on $x$ then both curves formed by
$\gamma$ and an arc in $x$ have about half the length of $x$, so they
have to be degenerate. If one bounds a M\"obius band then the core
would be a better choice for $x$, so necessarily in this case each
bounds a punctured disk, so the component of $\Sigma-x$ containing
$\gamma-x$ is necessarily a twice punctured disk. Next, if $\Gamma$
contains two short arcs $\gamma_1$ and $\gamma_2$ with each pair of
endpoints close along $x$, and the endpoints of $\gamma_1$ are not
nearly antipodal from the endpoints of $\gamma_2$, then we find $x'$
that bounds a twice punctured disk as in the first bullet.

There are now two cases that remain.

\item $\Gamma$ has one short arc on one side of
$x$ (which is a twice punctured disk) and the endpoints are far apart,
and on the other side it has two short arcs, with endpoints close, but
nearly antipodal from each other. Now use the arc from the twice
punctured disk side and an arc on the other side to construct $x'$
which is a better choice than $x$. 

\item Both sides of $x$ are twice punctured disks (and so $S$ is a
  4-times punctured sphere). If necessary, homotope $x$ to a
  geodesic. Now each side is a twice punctured disk with totally
  geodesic boundary. Double branched cover, branched over the
  punctures, is a flat annulus, and the punctured disk is the quotient
  of the obvious involution with two fixed points. The short arc of
  $\Gamma$ is, up to isotopy, any of the geodesic arcs perpendicular
  to both boundary components. There is a 1-parameter family of such
  arcs and we may arrange that on one side the endpoints are close,
  and on the other side at least one of the two endpoints is close to
  the first two. Now we may construct a scc whose length is only
  slightly longer than $\frac 12length(x)$ and bounds a twice punctured
  disk. 
\end{itemize}

{\bf Step 5.} Let $M$ be the maximal number of 
non-degenerate, pairwise
disjoint, nonparallel scc's in $\Sigma\setminus \{ \mbox{ punctures
}\}$. We claim that there is an essential annulus whose boundary
components are at least $\frac{\eta_0}M \sqrt{Area(\Sigma)}$ apart. 
To see this, consider a 1-Lipschitz map $\Sigma\to [0,d(x,y)]$ that
sends $x$ to 0 and $y$ to $d(x,y)$. Subdivide $[0,d(x,y)]$ into $M$
subintervals and consider the $M+1$ preimages of the vertices. Without
loss the map is transverse to the vertices. First note that each
preimage must contain a component which is non-degenerate
 (otherwise $x$
and $y$ would be separated by a union of degenerate curves). By the
choice of $M$, among them two are parallel.

{\bf Step 6.} (Besicovich Lemma) Let $A$ be an annulus and $f:A\to
[0,d]$ a 1-Lipschitz map. Then there is a core curve in $A$ of length
$\leq Area(A)/d$. Imagining that $f$ is smooth, this can be seen by
defining a 2-form $\omega$ on $A$ that evaluates on $(v,w)$ where
$df(v)=0$ as $\pm |v|df(w)$ (the sign is determined by
orientations). Then the 1-Lipschitz condition implies that $\omega$ is
a subarea form, and hence
$$Area(A)\geq \int_A\omega=\int_0^d length(f^{-1}(t))dt$$ which means
that for some $t$, $length(f^{-1}(t))\leq Area(A)/d$. 

Put $c=f^{-1}(t)$. Then 
$c$ is an essential scc with 
length $\le M \sqrt{Area(\Sigma)}/ \eta_0$ and 
width $ \ge \eta_0 \sqrt{Area(\Sigma)}/M$.
\end{proof}

\begin{remark}\label{modulus}
We could define modulus of a curve $c$ by taking the sup of
$\frac{I(c,x)}{l(x)}$ as $x$ runs over all {\it immersed} curves in
the surface. Call the resulting number $m'(c)$. Then {\it a priori}
$m(c)\le m'(c)$, but the proof of Lemma \ref{1} shows that $m'(c)\le
R/\sqrt{I(a,b)}$.
\end{remark}

In what follows it is best to think of 1-sided curves as being
interchangeable with the associated degenerate curves. Sometimes one
is more suitable than the other. To formalize this, consider the
complex $X'$ whose vertices are isotopy classes of essential 2-sided
scc's and simplices correspond to collections of curves that can be
realized disjointly. Then $X\cong X'$ by the map which is ``identity''
except it sends 1-sided curves to the associated degenerate 2-sided
curves. 

\begin{remark}\label{m2l}
Note that if $c,d$ are 2-sided, $m(c)\leq
R/\sqrt{I(a,b)}$ and $l(d)\leq R\sqrt{I(a,b)}$ then $c$ and $d$ are
close in $X'$. Indeed, $I(c,d)\leq l(d)m(c)\leq R^2$, so
$dist(c,d)\leq 1+I(c,d)\leq 1+R^2$. The same holds if say $l(d)\leq
100R\sqrt{I(a,b)}$. 
\end{remark}

The {\it coarse midpoint} for $\{a,b\}$ is the set $Mid(a,b)\subset
X'$ of all 2-sided essential curves $c$ with $m_{ab}(c)\leq
R/\sqrt{I(a,b)}$. 
In view of Remark \ref{m2l}, we could also work with the set
$Mid'(a,b)$ of all curves $c$ with $l_{ab}(c)\leq R\sqrt{I(a,b)}$. The
two sets intersect and any point in $Mid$ is within $1+R^2$ of any
point in $Mid'$. Thus both sets have diameter at most $2+2R^2$.

Recall that in the above discussion $a$ and $b$ are curve
systems. For example, note that $Mid(Na,Nb)=Mid(a,b)$. We now define
the {\it weighted midpoint} $Mid(a,b;\frac pq)$ as $Mid(qa,pb)$. Thus
$Mid(a,b)=Mid(a,b;1)$. Finally, we define the {\it coarse
  geodesic} between $a$ and $b$ as
$$\Lambda_{ab}=\bigcup\{Mid(a,b;\frac pq)|0< \frac pq <\infty\}$$
Note that for $N$ large any component of $a$ is in $Mid(a,b;\frac 1N)$
and any component of $b$ is in $Mid(a,b;N)$.

\begin{lemma}\label{fellow travel}
Suppose $I(b,c)=0$ and $I(a,b)=I(a,c)=I$. Then $Mid(a,b)$ and
$Mid(a,c)$ are close (here $a,b,c$ are 2-sided curve systems).
\end{lemma}

\begin{proof}
Let $x$ be a curve with $m_{ab}(x)\leq R/\sqrt I$ and $l_{ab}(x)\leq
R\sqrt I$. Then
$$I(c,x)\leq m_{ab}(x)l_{ab}(c)\leq \frac R{\sqrt I}I=R\sqrt I$$
so $l_{ac}(x)\leq 2R\sqrt I$ and the claim follows from Remark \ref{m2l}.
\end{proof}

Now suppose $a,b,c$ are three curve systems, each pair filling. Consider
the curve systems $$\tilde a=I(b,c)a,\quad\tilde b=I(c,a)b,\quad\tilde
c=I(a,b)c$$ and note that the intersection number between each pair is
$\tilde I=I(a,b)I(a,c)I(b,c)$.

\begin{lemma}\label{centers}
The three sets $Mid(\tilde a,\tilde b),Mid(\tilde b,\tilde
c),Mid(\tilde c,\tilde a)$ are close to each other.
\end{lemma}

\begin{proof}
Let $x\in Mid(\tilde a,\tilde b)$ be a curve from Lemma \ref{1}
(recall that if $x$ bounds a M\"obius band the core curve is in
$Mid(\tilde a,\tilde b)$). We
estimate the length of $x$ with respect to $\tilde a,\tilde c$.
$I(x,\tilde a)\leq R\sqrt{\tilde I}$ by the choice of $x$, and the
other estimate is as in Remark \ref{m2l}: 
$$I(x,\tilde c)\leq m_{\tilde a\tilde b}(x)l_{\tilde a\tilde b}(\tilde
c)\leq \frac R{\sqrt{\tilde I}}(I(\tilde a,\tilde c)+I(\tilde b,\tilde
c))=2R\sqrt {\tilde I}$$ Thus $$l_{\tilde a,\tilde c}(x)=I(x,\tilde
a)+I(x,\tilde c)\leq 3R\sqrt{\tilde I}$$ and we are done by Remark
\ref{m2l}. 
\end{proof}

By $Center(a,b,c)$ denote the union $Mid(\tilde a,\tilde b)\cup
Mid(\tilde b,\tilde c)\cup Mid(\tilde c,\tilde a)$. It is a nonempty
set of uniformly bounded diameter.

\begin{lemma}\label{thin triangles}
If $M\geq N$ then $Mid(M\tilde a,N\tilde b)$ is close to $Mid(M\tilde
a,N\tilde c)$.
\end{lemma}

\begin{proof}
Let $x$ be a special curve from Lemma \ref{1} for the pair $M\tilde
a,N\tilde b$. Then
$$I(x,N\tilde c)\leq m_{M\tilde a,N\tilde b}(x)l_{M\tilde a,N\tilde
  b}(N\tilde c)\leq \frac R{\sqrt{MN\tilde I}}(MN\tilde I+N^2\tilde
  I)\leq 2R\sqrt{MN\tilde I}$$ and we deduce as before that
$x$ is close to $Mid(M\tilde
a,N\tilde c)$.
\end{proof}

The following lemma is motivated by Axiom (3) (see below). If $x\in
Mid(a,b)$ we wish to show that $Center(a,b,x)$ and $x$ are close. The
issue is that $x$ may intersect one of the curves, say $a$, more often
than the other, so to compute $Center(a,b,x)$ we have to pass to
multiples of $a$ and $b$, with the coefficient of $a$ being smaller
than that of $b$. The point is that the effect on $x$ is that it gets
relatively shorter that way.

\begin{lemma}\label{ax3}
Suppose $l_{ab}(x)\leq R{\sqrt{I(a,b)}}$, i.e. $x\in Mid'(a,b;1)$. Then $x$ is
close to $Mid'(a,b;t)$ for all $t=\frac pq$ between $1$ and $\frac
{I(a,x)}{I(b,x)}$.
\end{lemma}

\begin{proof}
For concreteness assume $$1\leq \frac pq\leq \frac
{I(a,x)}{I(b,x)}$$ Then 
$l_{qa,pb}(x)=I(qa,x)+I(pb,x)=qI(a,x)+pI(b,x)$ and we wish to show 
$$qI(a,x)+pI(b,x)\leq CR\sqrt{I(qa,pb)}$$ for some constant $C>0$.
Squaring and writing $I=I(a,b)$ we need to show
$$q^2I(a,x)^2+p^2I(b,x)^2+2pqI(a,x)I(b,x)\leq C^2R^2pqI$$
Since $l_{ab}(x)\leq R\sqrt I$ we certainly have $I(a,x)\leq R\sqrt I$
and $I(b,x)\leq R\sqrt I$, so the last term is bounded (with $C^2=2$
on the right hand side). The first term can be bounded using $q\leq p$
and $I(a,x)\leq R\sqrt I$, and the second using $pI(b,x)\leq qI(a,x)$
and then estimating as in the last term.
\end{proof}

\subsection{Hyperbolicity criterion}

Let $X$ be a connected graph with the edge-path metric. Suppose that
for each pair $a,b$ of vertices we have a subset $\Lambda_{ab}$ with a
``coarse order'' (total, but not anti-symmetric), and
$\Lambda_{ab}=\Lambda_{ba}$ as sets but with reversed order. For
$x,y\in \Lambda_{ab}$ with $x\leq y$ write
$$\Lambda_{ab}[x,y]=\Lambda_{ab}[y,x]=\{z\in\Lambda_{ab}|x\leq z\leq
y\}$$ 
Also assume that $\phi:X^0\times X^0\times X^0\to X^0$ is a ternary
function on the vertices, symmetric under permutations, and such that
$\phi(a,b,c)\in \Lambda_{ab}$. We assume that
there is $K>0$ so that
\begin{enumerate}[(1)]
\item
  $HausDist(\Lambda_{ab}[a,\phi(a,b,c)],\Lambda_{ac}[a,\phi(a,b,c)])\leq K$
\item If $c,d$ are adjacent vertices then ${\rm{diam}}
  \Lambda_{ab}[\phi(a,b,c),\phi(a,b,d)]\leq K$
\item If $c\in\Lambda_{ab}$ then $\Lambda_{ab}[c,\phi(a,b,c)]\leq K$.
\end{enumerate}

\begin{thm}[Prop 3.1
\cite{bowditch:hyp}]
 It follows that $X$ is Gromov hyperbolic.
There exists $Q$ such that for all $a,b$,
$\Lambda_{ab}$
is in the $Q$-neighborhood of a geodesic
between $a,b$ (we say $\Lambda_{ab}$ is 
$Q$-quasi-convex).
\end{thm}

\subsection{Verifications}
We put $\phi=Center$. 
Axiom (1) follows from Lemma \ref{thin triangles}.  Lemma \ref{ax3}
reduces the verification of Axiom (3) to the case $I(a,c)=I(b,c)$ when
it follows immediately from definitions.

We now verify Axiom (2). For simplicity, we normalize so that
$I(a,b)=I(a,d)=I(b,d)=I$ and we assume that $I(a,c)\geq
I(b,c)$. Consider the point in $Center(a,b,c)$ that belongs to
$Mid(I(a,b)c,I(b,c)a)$ (such a point is suggestively called
$I(a,b)c+I(b,c)a$ in the figure below and it is represented by a
``bull's eye''). By Lemma \ref{fellow travel} there is a point in
$\Lambda[a,d]$ close to it, and a computation\footnote{Note that Lemma
  \ref{fellow travel} applies to $qa,pb,pc$, therefore $Mid(a,b;\frac
  pq)$ is uniformly close to $Mid(a,c,\frac pq)$ for any $p,q$. This
  gives a ``linear'' transfer map:
$I(a,b)c+I(b,c)a=Ic+I(b,c)a=I(a,c)\frac I{I(a,c)}c+I(b,c)a$ which
transfers to $I(a,c)d+I(b,c)a$} shows that this point is
$I(a,c)d+I(b,c)a$, using the same suggestive notation. Analogously, we
can start with the point $I(a,b)c+I(a,c)b$ also representing
$Center(a,b,c)$ and transfer it to $I(b,c)d+I(a,c)b$.

\centerline{\includegraphics{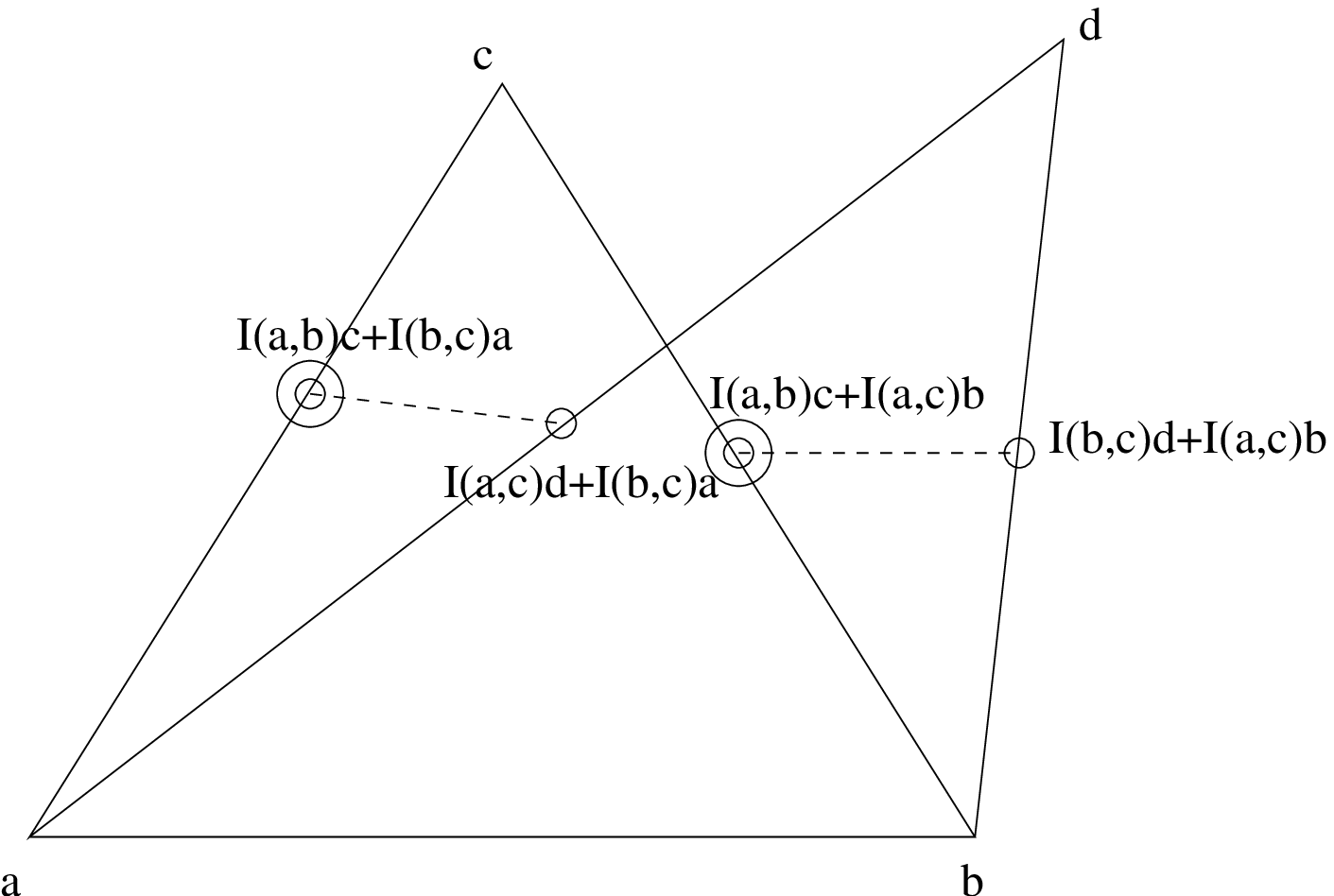}}

The two bull's eyes are close (they are in $Center(a,b,c)$), and
consequently the other two points are close. Note that one of them is
``too high'' and the other is ``too low'' (the center $Center(a,b,d)$
is represented by $a+d$ or $b+d$).

There is a 1-parameter family of points that interpolate between the
two points above. By $c_t$ denote the curve system
$tc\cup (1-t)d$, so $c_1=c$ and $c_0=d$. Note that
$I(c,c_t)=I(d,c_t)=0$ for all $t$ because $c,d$ are 2-sided. Take the point
$I(a,b)c+I(b,c)a$, transfer it to $\Lambda_{a,c_t}$, then to
$\Lambda_{b,c_t}$, and finally to $\Lambda_{bd}$. When $t=1$ we get
$I(b,c)d+I(a,c)b$ and when $t=0$ we get $I(a,c)d+I(b,c)b$, so for some
$0\leq t\leq 1$ we get the center for $abd$, proving the claim. The
only thing we have to ensure is that the transfer from
$\Lambda_{a,c_t}$ to $\Lambda_{b,c_t}$ is possible, i.e. that the
point in question is on the ``$c_t$-side'' of the centerpoint for
$abc_t$. A little calculation\footnote{
In a triangle $xyz$ the point $\alpha x+\beta
y=\frac{\alpha}{I(y,z)}I(y,z)x+
\frac{\beta}{I(x,z)}I(x,z)y$ is on the $x$-side of the center iff
$\frac{\alpha}{I(y,z)}\geq\frac{\beta}{I(x,z)}$ iff $\alpha I(x,z)\geq
\beta I(y,z)$. Applying this to $c_tab$ we get the condition
$\frac{I\times I(a,c)}{tI(a,c)+(1-t)I}I(c_t,b)\geq I(b,c)I(a,b)$ i.e. 
$I\times I(a,c)(tI(b,c)+(1-t)I)\geq I(b,c)I(tI(a,c)+(1-t)I)$ which
simplifies to $I(a,c)\geq I(b,c)$}
shows that the condition is that 
$I(a,c)\geq I(b,c)$.

\subsection{The action of $MCG(S)$ on $X$}
Let $MCG(S)$ be the mapping class group of $S$.  Then $MCG(S)$ acts on
$X$ by isometries.  Masur-Minsky (Prop 4.6 \cite{MM}) showed a
pseudo-Anosov element $g \in MCG(S)$ acts as a hyperbolic isometry on
$X$ in the following sense: there exists a quasi-geodesic $\gamma$ in
$X$ which is invariant by $g$. In this case $\gamma$ is called a
(quasi-geodesic) axis of $g$.  It is trivial that $g$ is not
hyperbolic if it is not pseudo-Anosov.

It is shown (Prop 11 of \cite{bf}) that the action of $MCG(S)$ on $X$
is weakly properly discontinuous (WPD) in the following 
sense: 
for every pseudo-Anosov element $g \in MCG(S)$, every $x \in X$, and 
every  $C \ge 0$, there exists $N >0$ such that 
the following set is finite:
$$\{h \in MCG(S)| dist(x,h(x)) \le C, dist(g^N(x),hg^N(x)) \le C \}.$$

We claim those results extend to non-orientable surfaces.
\begin{thm}
Let $\Sigma$ be a compact surface with $n$ punctures.
We exclude a sphere with $n \le 4$, a torus with $n \le 1$,
a projective plane with $n \le 2$ and a Klein bottle with 
$n \le 1$.
Then an element $g \in MCG(S)$ acts as a hyperbolic
isometry on the curve complex $X$ if and only if 
$g$ is pseudo-Anosov.
The action of $MCG(S)$ on $X$ is WPD.
\end{thm}

\begin{proof}
For the first part, the argument in \cite{MM}
uses the theory of measured lamination on $S$.
The theory does not require orientability of $S$, \cite{thurston}.

For the second part, the argument in \cite{bf}
also uses measured laminations, and does not assume 
orientability of $S$.
\end{proof}

We remark that Bowditch \cite{bowditch:tight} showed that the action 
of $MCG(S)$ on $X$ is ``acylindrical", which implies WPD.
His argument uses ``tight geodesics" on $X$, which is 
developed in \cite{minsky-masur:cc2} in the orientable case.

\end{document}